\documentclass[10pt,reqno]{amsart}
\usepackage{amssymb,amscd,amsbsy}
\usepackage{amsthm}
\renewcommand{\a}{\alpha}

\newcommand{\z}{\zeta}

\newcommand{\s}{\sigma}
\renewcommand{\t}{\tau}

\newcommand{\f}{\varphi}
\renewcommand{\o}{\omega}

\newcommand{\D}{\Delta}
\renewcommand{\L}{\Lambda}

\newcommand{\T}{{\Bbb T}}

\newcommand{\R}{{\Bbb R}}

\newcommand{\bS}{{\boldsymbol S}}

\newcommand{\rf}[1]{(\ref{#1})}

\newcommand{\df}{\stackrel{\mathrm{def}}{=}}

\newcommand{\const}{\operatorname{const}}

\newcommand{\eeq}{\end{equation}}
\newcommand{\beq}{\begin{equation}}
\newcommand{\bay}{\begin{eqnarray}}
\newcommand{\ba}{\begin{align*}}
\newcommand{\ea}{\end{align*}}
\newcommand{\ey}{\end{eqnarray}}
\newcommand{\bey}{\begin{eqnarray*}}
\newcommand{\eey}{\end{eqnarray*}}

\newcommand{\be}{\infty}

\newtheorem{thm}{\hspace{\parindent}Theorem}[section]

\theoremstyle{remark}

\newtheorem*{rem*}{Remark}

\begin{document}

\newcommand{\vse}{\vspace{.2in}}
\numberwithin{equation}{section}

\title{Functions  of perturbed operators}
\author{Aleksei Aleksandrov and Vladimir Peller}


\maketitle

\newcommand{\mt}{{\mathcal T}}

\footnotesize

{\bf Abstract.} We prove that if $0<\a<1$ and $f$ is in the H\"older class $\L_\a(\R)$, then for arbitrary self-adjoint operators $A$ and $B$
with bounded $A-B$, the operator $f(A)-f(B)$ is bounded and $\|f(A)-f(B)\|\le\const\|A-B\|^\a$. We prove a similar result for functions 
$f$ of the Zygmund class $\L_1(\R)$: $\|f(A+K)-2f(A)+f(A-K)\|\le\const\|K\|$, where $A$ and $K$ are self-adjoint operators. Similar results also hold for 
all H\"older-Zygmund classes $\L_\a(\R)$, $\a>0$. We also study properties of the operators $f(A)-f(B)$ for $f\in\L_\a(\R)$ and self-adjoint operators $A$ and $B$ such that $A-B$ belongs to the Schatten--von Neumann class $\bS_p$. We consider the same problem for higher order differences.
Similar results also hold for unitary operators and for contractions.

\medskip

\begin{center}
{\bf\large Fonctions d'op\'erateurs perturb\'es}
\end{center}

\medskip

{\bf R\'esum\'e.} Nous montrons que si $0<\a<1$ et $f$ appartient \`a la classe de H\"older $\L_\a(\R)$, alors pour tous les op\'erateurs 
$A$ et $B$ auto-adjoints dont la diff\'erence est born\'ee on a: $\|f(A)-f(B)\|\le\const\|A-B\|^\a$. Nous obtenons un r\'esultat similaire pour 
les fonctions de la classe de Zygmund $\L_1(\R)$: $\|f(A+K)-2f(A)+f(A-K)\|\le\const\|K\|$, o\`u $A$ et $K$ sont des op\'erateurs auto-adjoints.
Un r\'esultat similaire est aussi vrai pour toutes les classes de H\"older--Zygmund $\L_\a(\R)$, $\a>0$. Nous \'etudions aussi les propri\'et\'es
des op\'erateurs $f(A)-f(B)$ si $f\in\L_\a(\R)$ et $A$ et $B$ sont des op\'erateurs auto-adjoints dont la diff\'erence appartient \`a la classe
de Schatten--von Neumann $\bS_p$. Nous consid\'erons le m\^eme probl\`eme pour les diff\'erences d'ordre arbitraire. On peut
obtenir des r\'esultats similaires pour les op\'erateurs unitaires et pour les contractions.

\normalsize

\

\begin{center}
{\bf\large Version fran\c caise abr\'eg\'ee}
\end{center}

\medskip

Il est bien connu qu'il y a des fonctions $f$ lipschitziennes sur la droite r\'eelle $\R$ qui ne sont pas lipschitziennes op\'eratorielles,
c'est-\`a-dire la condition
$$
|f(x)-f(y)|\le\const|x-y|,\quad x,~y\in\R,
$$
n'implique pas que pour tous les op\'erateurs auto-adjoints $A$ et $B$ l'in\'egalit\'e
$$
\|f(A)-f(B)\|\le\const\|A-B\|.
$$
soit vraie.

Il se trouve que la situation change dramatiquement si l'on consid\`ere les fonctions de la classe $\L_\a(\R)$ de H\"older d'ordre $\a$, $0<\a<1$.
Nous montrons que si $A$ et $B$ sont des op\'erateurs auto-adjoints dans un espace hilbertien et $f\in\L_\a(\R)$
(c'est-\`a-dire $|f(x)-f(y)|\le\const|x-y|^\a$), alors
$$
\|f(A)-f(B)\|\le\const\|A-B\|^\a.
$$

On peut consid\'erer un probl\`eme similaire pour la classe de Zygmund $\L_1(\R)$ de fonctions continues sur $\R$ telles que
$$
|f(x+t)-2f(x)+f(x-t)|\le\const|t|,\quad x,\,t\in\R.
$$
Nous \'etablissons que dans ce cas $f$ doit satisfaire \`a l'in\'egalit\'e
$$
\|f(A+K)-2f(A)+f(A-K)\|\le\const\|K\|
$$
o\`u $A$ et $K$ sont des op\'erateurs auto-adjoints. 

Nous consid\'erons aussi les espaces $\L_\a(\R)$ pour tous les $\a>0$. Si $\a>0$, la classe $\L_\a(\R)$ consiste en fonctions $f$ continues telles que
$$
\left|\sum_{k=0}^n(-1)^{n-k}\left(\begin{matrix}n\\k\end{matrix}\right)f(x+kt)\right|\le\const|t|^\a\quad(\mbox{ici}\quad n-1\le\a <n).
$$
Nous montrons que dans ce cas
$$
\left\|\sum_{j=0}^n(-1)^{n-j}\left(\begin{matrix}n\\j\end{matrix}\right)f\big(A+jK\big)\right\|\le\const\|K\|^\a,\quad n>\a,
$$
pour chaque op\'erateurs $A$ et $K$ auto-adjoints.

De plus, on peut g\'en\'eraliser les r\'esultats ci-dessus pour les op\'erateurs $A$ pas n\'ecessairement
born\'es (le th\'eor\`eme \ref{vp} dans la version anglaise).

On peut consid\'erer les m\^emes probl\`emes pour les op\'erateurs unitaires et pour les contractions. 

Nous consid\'erons aussi le cas de modules de continuit\'e arbitraires. Soit $\o$ une fonction croissante continue sur $[0,\be)$ telle que
$\o(0)=0$ et $\o(x+y)\le\o(x)+\o(y)$, $x,\,y\ge0$. Nous d\'efinissons la fonction $\o^*$ par
$$
\o^*(x)=x\int_x^\be\frac{\o(t)}{t^2}\,dt,\quad x\ge0.
$$
Notons par $\L_\o(\R)$ l'espace de fonctions $f$ sur $\R$ telles que
$$
|f(x)-f(y)|\le\const\o(|x-y|),\quad x,y\in\R.
$$
Nous d\'emontrons que si $A$ et $B$ sont des op\'erateurs auto-adjoints dont la diff\'erence $A-B$ est born\'ee et $f\in\L_\o$, alors
$$
\|f(A)-f(B)\|\le\const\o^*(\|A-B\|).
$$
On peut obtenir un r\'esultat similaire pour les modules de continuit\'e d'ordre arbitraire. Nous avons aussi obtenu les m\^emes r\'esultats pour les 
op\'erateurs unitaires et pour les contractions.

Maintenant nous consid\'erons le probl\`eme suivant. Rappelons que $\bS_p$ est la classe de Schatten--von Neumann constitu\'ee des
op\'erateurs $T$ dans un espace hilbertien pour lesquels les nombres singuliers $s_n(T)$ appartiennent \`a l'espace $\ell^p$.

Supposons que $f\in\L_\a(\R)$, $0<\a<1$, et $p>1$. Nous montrons que pour chaque op\'erateurs auto-adjoints $A$ et $B$ dont la diff\'erence $A-B$ appartient \`a  $\bS_p$, on a
$$
f(A)-f(B)\in\bS_{\frac p\a}\quad\mbox{et}\quad\|f(A)-f(B)\|_{\bS_\frac p\a}\le\const\|A-B\|^\a.
$$
Si $p=1$ et $f\in\L_\a(\R)$, $0<\a<1$, la condition $A-B\in\bS_1$ implique que
$$
f(A)-f(B)\in\bS_{\frac 1\a,\be},
$$
o\`u l'espace $\bS_{q,\be}$ est form\'e des op\'erateurs dont les nombres singuliers satisfont \`a la condition
$$
\sup_{n\ge0}(1+n)^{1/q}s_n(T)<\be.
$$
Nous avons aussi obtenu des analogues des r\'esultats ci-dessus pour les diff\'erences d'ordre arbitrare.
On peut obtenir des r\'esultats similaires pour les op\'erateurs unitaires et pour
les contractions.

\begin{center}
------------------------------
\end{center}

\setcounter{section}{0}
\section{\bf Introduction}
\setcounter{equation}{0}

\medskip

It is well known that a Lipschitz function on the real line is not necessarily {\it operator Lipschitz}, i.e.,
the condition 
$$
|f(x)-f(y)|\le\const|x-y|,\quad x,~y\in\R,
$$
does not imply that for self-adjoint operators
$A$ and $B$ on Hilbert space,
$$
\|f(A)-f(B)\|\le\const\|A-B\|.
$$
The existence of such functions was proved in \cite{F} (see also \cite{F2} and \cite{K}). 
Later in \cite{Pe1} necessary conditions were found for a function $f$ to be operator Lipschitz. Those necessary conditions imply that Lipschitz  functions do not have to be operator Lipschitz.
It is also well known that a continuously differentiable function does not have to be operator differentiable, see \cite{Pe1} and \cite{Pe2}.
Note that the necessary conditions obtained in \cite{Pe1} and \cite{Pe2} are based on the nuclearity criterion for Hankel operators, see \cite{Pe3}.

It turns out that the situation dramatically changes if we consider H\"older classes $\L_\a(\R)$ with 
\mbox{$0<\a<1$}. In this case such functions are necessarily {\it operator H\"older of order $\a$}, i.e., the condition
$$
|f(x)-f(y)|\le\const|x-y|^\a,\quad x,~y\in\R,
$$
implies that for self-adjoint operators
$A$ and $B$ on Hilbert space,
\bay
\label{fo}
\|f(A)-f(B)\|\le\const\|A-B\|^\a.
\ey
Moreover, a similar result holds for the Zygmund class $\L_1(\R)$, i.e., the fact that 
$$
|f(x+t)-2f(x)+f(x-t)|\le\const|t|,\quad x,\,t\in\R,
$$
and $f$ is continuous implies that $f$ is {\it operator Zygmund}, i.e., for self-adjoint operators $A$ and $K$,
\bay
\label{so}
\|f(A+K)-2f(A)+f(A-K)\|\le\const\|K\|.
\ey
We also obtain similar results for the whole scale of H\"older--Zygmund classes $\L_\a(\R)$ for $0<\a<\be$. Recall that for $\a>1$,
the class $\L_\a(\R)$ consists of continuous functions $f$ such that
$$
\left|\sum_{k=0}^n(-1)^{n-k}\left(\begin{matrix}n\\k\end{matrix}\right)f(x+kt)\right|\le\const|t|^\a,\quad\mbox{where}\quad n-1\le\a <n.
$$

The same problems can be considered for unitary operators and for functions on the unit circle, and for contractions and 
analytic functions in the unit disk.
 
 To prove \rf{fo}, we use
a crucial estimate  obtained for trigonometric polynomials and unitary operators in \cite{Pe1}  and for entire functions of exponential type and self-adjoint operators in \cite{Pe2}. We state here the result for self-adjoint operators. It can be considered as an analog of Bernstein's inequality.

{\it Let $f$ be an entire function of exponential type $\s$ that is bounded on the real line $\R$. Then for self-adjoint operators $A$ and $B$ with bounded $A-B$ the following inequality holds:}
\bay
\label{nB}
\|f(A)-f(B)\|\le\const\s\|f\|_{L^\be(\R)}\|A-B\|.
\ey

Inequality \rf{nB} was proved by using double operator integrals and the Birman--Solomyak formula:
$$
f(A)-f(B)=\iint\frac{f(x)-f(y)}{x-y}\,dE_A(x)(A-B)\,dE_B(y),
$$
where $E_A$ and $E_B$ are the spectral measures of self-adjoint operators $A$ and $B$; we refer the reader to \cite{BS1}, \cite{BS2} and \cite{BS3} for the theory of double operator integrals. Note that $A$ and $B$ do not have to be bounded, but $A-B$ must be bounded.

To estimate the second difference \rf{so}, we use the corresponding analog
of Bernstein's inequality which was obtained in \cite{Pe4} with the help of
triple operator integrals.
 To estimate higher order differences, we need multiple operator integrals. We refer the reader to 
\cite{Pe4} for definitions and basic results on multiple operator integrals.

We also consider in this paper the problem of the behavior of functions of operators $f(A)$ under perturbations of $A$ by operators of Schatten--von Neumann class $\bS_p$ in the case when $f\in\L_\a(\R)$. 

\medskip

\section{\bf Norm estimates for unitary operators}
\setcounter{equation}{0}

\medskip

We start with first order differences. We use the notation by $\L_\a$, $0<\a<\be$, for the scale of H\"older--Zygmund classes on the unit circle $\T$.

\begin{thm}
\label{uH}
Let $0<\a<1$. Then there is a constant $c>0$ such that for 
every $f\in\L_\a$ and for arbitrary unitary operators $U$ and $V$ on Hilbert space the following inequality holds:
$$
\|f(U)-f(V)\|\le c\|f\|_{\L_\a}\cdot\|U-V\|^\a.
$$
\end{thm}

\begin{thm}
\label{oLu}
There exists a constant $c>0$ such that for every function $f\in\L_1$ and for arbitrary unitary operators $U$ and $V$ on Hilbert space
the following inequality holds:
$$
\|f(U)-f(V)\|\le c\|f\|_{\L_1}\left(2+\log_2\frac1{\|U-V\|}\right)\|U-V\|.
$$
\end{thm}

Note that this result improves an estimate  obtained in \cite{F} for Lipschitz functions in the case of bounded self-adjoint operators.

We proceed now to higher order differences.


\begin{thm}
\label{hou}
Let $n$ be a positive integer and $0<\a<n$. Then there exists a constant $c>0$ such that for every $f\in\L_\a$ and
for an arbitrary unitary operator $U$ and an arbitrary bounded self-adjoint operator $A$ on Hilbert space the following inequality holds:
$$
\left\|\sum_{k=0}^n(-1)^{n-k}\left(\begin{matrix}n\\k\end{matrix}\right)f\big(e^{{\rm i}kA}U\big)
\right\|\le c\|f\|_{\L_\a}\|A\|^\a.
$$
\end{thm}

Let us consider now a more general problem. Suppose that $\o$ is a modulus of continuity, i.e., $\o$ is a nondecreasing continuous function on $[0,\be)$
such that $\o(0)=0$ and $\o(x+y)\le\o(x)+\o(y)$, $x,\,y\ge0$. The space $\L_\o$ consists of functions $f$ on $\T$ such that
$$
|f(\z)-f(\t)|\le\const\o(|\z-\t|),\quad\z,\,\t\in\T.
$$
With a modulus of continuity $\o$ we associate the function $\o^*$ defined by
$$
\o^*(x)=x\int_x^\be\frac{\o(t)}{t^2}\,dt,\quad x\ge0.
$$

\begin{thm}
\label{om}
Suppose that $\o$ is a modulus of continuity and $f\in\L_\o$. If $U$ and $V$ are unitary operators, then
$$
\|f(U)-f(V)\|\le\const\|f\|_{\L_\o}\,\o^*(\|U-V\|).
$$
\end{thm}

In particular, if $\o^*(x)\le\const\o(x)$, then for unitary operators $U$ and $V$
$$
\|f(U)-f(V)\|\le\const\|f\|_{\L_\o}\,\o(\|U-V\|).
$$

We have also proved an analog of Theorem \ref{om} for higher order differences.

\medskip

\section{\bf Norm estimates for contractions}
\setcounter{equation}{0}

\medskip

We denote here by $(\L_\a)_+$ the set of functions $f\in\L_\a$, for which the Fourier coefficients $\hat f(n)$ vanish for $n<0$.

Recall that an operator $T$ on Hilbert space is called a contraction if $\|T\|\le1$. The following result is an analog of Theorem \ref{hou}
for contractions.

\begin{thm}
\label{conh}
Let $n$ be a positive integer and $0<\a<n$. Then there exists a constant $c>0$ such that for every $f\in(\L_\a)_+$ and
for arbitrary contractions $T$ and $R$ on Hilbert space, the following inequality holds:
$$
\left\|\sum_{k=0}^n(-1)^{n-k}\left(\begin{matrix}n\\k\end{matrix}\right)f\left(T+\frac{k}{n}(T-R)\right)
\right\|\le c\|f\|_{\L_\a}\|T-R\|^\a.
$$
\end{thm}



Note that an analog of Theorem \ref{om} also holds for contractions.

\medskip

\section{\bf Norm estimates for self-adjoint operators}
\setcounter{equation}{0}

\medskip

\begin{thm}
\label{pp}
Let $0<\a<1$ and let $f\in\L_\a(\R)$. Suppose that $A$ and $B$ are self-adjoint operators such that $A-B$ is bounded. Then $f(A)-f(B)$ is bounded and
$$
\|f(A)-f(B)\|\le\const\|f\|_{\L_\a(\R)}\|A-B\|^\a.
$$
\end{thm}

In this connection we mention the paper
\cite{F} where it was proved that for self-adjoint operators $A$ and $B$ with spectra in an interval $[a,b]$ and a function $\f\in\L_\a(\R)$,
the following inequality holds:
$$
\|\f(A)-\f(B)\|\le\const\|\f\|_{\L_\a(\R)}\left(\log\left(\frac{b-a}{\|A-B\|}+1\right)+1\right)^2\|A-B\|^\a
$$
(see also \cite{FN} where the above inequality is generalized for general moduli of continuity).

\begin{thm}
\label{vp}
Suppose that $n$ is a positive integer and  $0<\a<n$. Let $A$ be a self-adjoint operator and let $K$ be a bounded self-adjoint operator.
Then the map
\bay
\label{ue}
f\mapsto\big(\D_K^nf\big)(A)\df\sum_{j=0}^n(-1)^{n-j}\left(\begin{matrix}n\\j\end{matrix}\right)f\big(A+jK\big)
\ey
has a unique extension from $L^\be\cap\L_\a(\R)$ to a  sequentially continuous operator from
$\L_\a(\R)$ (equipped with the weak-star topology) to the space of bounded linear operators on Hilbert space (equipped with the strong operator topology) and
$$
\|\big(\D_K^nf\big)(A)\|\le\const\|f\|_{\L_\a(\R)}\|K\|^\a.
$$
\end{thm}

We use the same notation $(\D_K^nf\big)(A)$ for the unique extension of the map \rf{ue}.

We can also prove an analog of Theorem \ref{om} for self-adjoint operators.

\medskip

\section{\bf Perturbations of class $\bS_p$}
\setcounter{equation}{0}

\medskip

In this section we consider the behavior of functions of self-adjoint operators under perturbations of Schatten--von Neumann class $\bS_p$.
Similar results also hold for unitary operators and for contractions.

Recall that the spaces $\bS_p$ and $\bS_{p,\be}$ consist of operators $T$ on Hilbert space such that 
$$
\|T\|_{\bS_p}\df\left(\sum_{n\ge0}\big(s_n(T)\big)^p\right)^{1/p}<\be\quad\mbox{and}\quad\|T\|_{\bS_{p,\be}}\df\sup_{n\ge0}(1+n)^{1/p}s_n(T)<\be.
$$

\begin{thm}
\label{pr}
Let $1\le p<\be$, $0<\a<1$, and let $f\in\L_\a(\R)$. Suppose that $A$ and $B$ are self-adjoint operators such that $A-B\in\bS_p$. Then
$$
f(A)-f(B)\in\bS_{\frac{p}{\a},\be}\quad\mbox{and}\quad
\|f(A)-f(B)\|_{\bS_{\frac{p}{\a},\be}}\le\const\|f\|_{\L_\a(\R)}\|A-B\|_{\bS_p}^\a.
$$
\end{thm}

Note that in Theorem \ref{pr} in the case $p>1$ we can replace the condition $A-B\in\bS_p$ with the condition $A-B\in\bS_{p,\be}$.

Using interpolation arguments, we can deduce from Theorem \ref{pr} the following result.

\begin{thm}
\label{pr1}
Let $1<p<\be$, $0<\a<1$, and let $f\in\L_\a(\R)$. Suppose that $A$ and $B$ are self-adjoint operators such that $A-B\in\bS_p$. Then
$$
f(A)-f(B)\in\bS_{\frac{p}{\a}}\quad\mbox{and}\quad
\|f(A)-f(B)\|_{\bS_{\frac{p}{\a}}}\le\const\|f\|_{\L_\a(\R)}\|A-B\|_{\bS_p}^\a.
$$
\end{thm}

Let us now state similar results for higher order differences.

\begin{thm}
\label{hod}
Suppose that $n$ is a positive integer, $\a$ is a positive number such that $n-1\le\a<n$, $f\in\L_\a(\R)$, and $n\le p<\be$. Let $A$ be a self-adjoint operator and let $K$ be a self-adjoint operator of class $\bS_p$.
Then the operator $\big(\D_K^nf\big)(A)$ defined in Theorem {\em\ref{vp}} belongs to $\bS_{\frac{p}{\a},\be}$
and
$$
\|\big(\D_K^nf\big)(A)\|_{\bS_{\frac{p}{\a},\be}}\le\const\|f\|_{\L_\a(\R)}\|K\|_{\bS_p}^\a.
$$
\end{thm}

\begin{thm}
\label{hodn}
Suppose that $n$ is a positive integer, $\a$ is a positive number such that $n-1\le\a<n$, $f\in\L_\a(\R)$, and $n<p<\be$. Let $A$ be a self-adjoint operator and let $K$ be a self-adjoint operator of class $\bS_p$.
Then the operator $\big(\D_K^nf\big)(A)$ defined in Theorem {\em\ref{vp}} belongs to $\bS_{\frac{p}{\a}}$
and
$$
\|\big(\D_K^nf\big)(A)\|_{\bS_{\frac{p}{\a}}}\le\const\|f\|_{\L_\a(\R)}\|K\|_{\bS_p}^\a.
$$
\end{thm}


\footnotesize

\noindent
\begin{tabular}{p{10cm}p{15cm}}
A.B. Aleksandrov & V.V. Peller \\
St-Petersburg Branch & Department of Mathematics \\
Steklov Institute of Mathematics  & Michigan State University \\
Fontanka 27, 191023 St-Petersburg & East Lansing, Michigan 48824\\
Russia&USA
\end{tabular}

\end{document}